\documentclass{amsart}
\usepackage{amssymb}
\def\version{Version W07-v3f last changed 23 August 2004 by PBG Printed \number
\day-\number\month-\number\year}
\date{\version}
\usepackage{a4}
\begin{document}
\newtheorem{theorem}{Theorem}[section]
\newtheorem{definition}[theorem]{Definition}
\newtheorem{lemma}[theorem]{Lemma}
\newtheorem{example}[theorem]{Example}
\newtheorem{remark}[theorem]{Remark}
\newtheorem{corollary}[theorem]{Corollary}
\def\ffrac#1#2{{\textstyle\frac{#1}{#2}}}
\def\qedbox{\hbox{$\rlap{$\sqcap$}\sqcup$}}
\font\pbglie=eufm10\def\gg{\text{\pbglie g}}
\def\id{\text{Id}}
\makeatletter
 \renewcommand{\theequation}{%
 \thesection.\alph{equation}}
 \@addtoreset{equation}{section}
 \makeatother
\title[Curvature homogeneous manifolds]
{Curvature homogeneous signature $(2,2)$ manifolds}
\author{C. Dunn, P. Gilkey and S. Nik\v cevi\'c}
\begin{address}{CD \& PG: Mathematics Department, University of Oregon,
Eugene Or 97403 USA.\newline Email: {\it cdunn@darkwing.uoregon.edu, gilkey@darkwing.uoregon.edu}}
\end{address}
\begin{address}{SN: Mathematical Institute, Sanu,
Knez Mihailova 35, p.p. 367,
11001 Belgrade,
Serbia and Montenegro.
\newline Email: {\it stanan@mi.sanu.ac.yu}}\end{address}
\begin{abstract}We exhibit a family of generalized plane wave manifolds of signature $(2,2)$. The geodesics in these manifolds
extend for infinite time (i.e. they are complete), they are spacelike and timelike Jordan Osserman, and they are spacelike and
timelike Jordan Ivanov-Petrova. Some are irreducible symmetric spaces. Some are homogeneous spaces but not symmetric. Some are
$1$-curvature homogeneous but not homogeneous. All are
$0$-modeled on the same irreducible symmetric space. We determine the Killing vector fields for these manifolds.
\end{abstract}
\maketitle
\section{Introduction}

We begin by recalling some definitions. Let $\mathcal{M}:=(M,g_M)$ be a pseudo-Riemannian manifold of signature $(p,q)$. 
Let $\nabla^kR_{\mathcal{M}}$ be the $k^{\operatorname{th}}$ covariant derivative of the curvature tensor and let
$G_{\mathcal{M}}$ be the isometry group. The manifold $\mathcal{M}$ is {\it locally symmetric} if $\nabla R_{\mathcal M}=0$ and the
manifold
$\mathcal{M}$ is {\it homogeneous} if $G_{\mathcal{M}}$ acts transitively on $\mathcal{M}$; note that a simply connected
complete local symmetric space is homogeneous. One says that
$\mathcal{M}$ is {\it spacelike} (resp. {\it timelike}) {\it Jordan Osserman} if the Jacobi operator has constant
Jordan normal form on the pseudo-sphere bundles of unit spacelike (resp. timelike) tangent vectors. Similarly, one says that a
pseudo-Riemannian manifold
$\mathcal{M}$ is {\it spacelike} (resp. {\it timelike}) {\it Jordan Ivanov-Petrova} if the skew-symmetric curvature operator has
constant Jordan normal form on the Grassmann bundles of oriented spacelike (resp. timelike) $2$ planes.

Examples are at the heart of modern Differential Geometry. Let $(x,y,\tilde x,\tilde y)$ be coordinates on $\mathbb{R}^4$. To
simplify the discussion, we shall only give the non-zero entries up to the usual symmetries when defining certain tensors.  Let
$\mathcal{A}(\mathcal{O})$ be the set of real analytic functions on a connected open subset
$\mathcal{O}\subset\mathbb{R}$. If
$f\in\mathcal{A}(\mathbb{R})$, define a pseudo-Riemannian metric of neutral signature $(2,2)$ on $\mathbb{R}^4$ by setting:
\begin{eqnarray*}
g_f(\partial_x,\partial_x)=-2f(y),\quad g_f(\partial_x,\partial_{\tilde x})=g_f(\partial_y,\partial_{\tilde y})=1\,.
\end{eqnarray*}
The {\it generalized plane wave} manifolds $\mathcal{M}_f:=(\mathbb{R}^4,g_f)$ are complete and form a rich family of examples;
we shall explore their geometry in some detail in this paper. These metrics are a special case of {\it Walker metrics}; see
\cite{CGM04} for related work.

Derdzinski \cite{D00} studied the manifolds $\mathcal{M}_f$ previously and showed:

\begin{theorem}[Derdzinski]\label{thm-1.1} 
Adopt the notation established above. Then:
\begin{enumerate}
\item If $f^{(3)}$ never vanishes, then $\alpha_2:=f^{(4)}f^{(2)}(f^{(3)})^{-2}$is an isometry invariant.
\item $\mathcal{M}_f$ is symmetric if and only if $f^{(3)}=0$.
\item $\mathcal{M}_f$ is curvature homogeneous if and only if $f^{(2)}=0$ identically or if $f^{(2)}$ never vanishes.
\end{enumerate}\end{theorem}
This result enabled him to conclude:
\begin{corollary}[Derdzinski]\label{cor-1.2} 
There exist neutral signature Ricci-flat $4$ dimensional pseudo-Riemannian manifolds
which are curvature-homogeneous but which are not locally homogeneous.\end{corollary}

In this paper, we will study these manifolds in further detail.

\begin{theorem}\label{thm-1.3}
 We have that:
\begin{enumerate}
\item The non-zero components of $\nabla^kR_{\mathcal{M}_f}$ are:
$$\nabla^kR_{\mathcal{M}_f}(\partial_x,\partial_y,\partial_y,\partial_x;\partial_y,...,\partial_y)=f^{(k+2)}\,.$$
\item All geodesics in $\mathcal{M}_f$ extend for infinite time.
\item If $P\in\mathbb{R}^4$, then $\exp_{\mathcal{M}_f,P}:T_P\mathbb{R}^4\rightarrow\mathbb{R}^4$ is
a diffeomorphism.
\item All the local scalar Weyl invariants of $\mathcal{M}_f$ vanish.
\item  If $f^{(2)}\ne0$, then $\mathcal{M}_f$ is spacelike and timelike Jordan Osserman.
\item If $f^{(2)}\ne0$, then $\mathcal{M}_f$ is spacelike and timelike Jordan
Ivanov-Petrova.
\item $\mathcal{M}_f$ is realizable as a hypersurface in $\mathbb{R}^{(2,3)}$.
\end{enumerate}
\end{theorem}

We extend the earlier result of Derdzinski:

\begin{theorem}\label{thm-1.4} We have that:
\begin{enumerate}\item $\mathcal{M}_f$ is symmetric if and only if $f^{(3)}\equiv0$.
\item $\mathcal{M}_f$ is homogeneous if and only if $f^{(2)}=ae^{\lambda y}$ for some $a,\lambda\in\mathbb{R}$.
\end{enumerate}\end{theorem}

\begin{remark}\label{rmk-1.5}
\rm Setting $f=e^y+e^{2y}$ yields a complete spacelike and timelike Jordan Osserman manifold of
signature $(2,2)$ which is not homogeneous. This is not possible in the Riemannian setting as Chi \cite{refChi} showed that a
complete $4$ dimensional Riemannian Osserman manifold is necessarily either a rank
$1$-symmetric space or is flat. The manifold $\mathcal{M}_{e^y+e^{2y}}$ is also a complete spacelike
and timelike Jordan Ivanov-Petrova manifold. Again, this is not possible in the Riemannian setting as Ivanov and Petrova
\cite{IvPe98} showed that a complete $4$ dimensional Riemannian Ivanov-Petrova manifold either has constant sectional curvature or is
a warped product of an interval with a manifold of constant sectional curvature. All the Weyl scalar invariants of
$\mathcal{M}_{e^y+e^{2y}}$ vanish. This is not possible in the Riemannian setting as Pr\"ufer, Tricerri, and Vanhecke \cite{PTV96}
showed that if all local scalar Weyl invariants up to order $\frac12m(m-1)$ are constant on a Riemannian manifold $\mathcal{M}$, then
$\mathcal{M}$ is locally homogeneous and
$\mathcal{M}$ is determined up to local isometry by these invariants. We note that there exist Lorentzian
manifolds all of whose Weyl scalar invariants vanish, see for example the discussion in Koutras and
McIntosh \cite{KM96} or Pravda, Pravdov\'a, Coley, and Milson \cite{PPCM02}.
\end{remark}

Let $\langle\cdot,\cdot\rangle$ be a non-degenerate inner product on a finite dimensional vector space $V$. Assume given
$A^i\in\otimes^{4+i}V^*$ for $i=0,1,...,k$. Set
$$\mathcal{U}^k:=(V,\langle\cdot,\cdot\rangle,A^0,A^1,...,A^k)\,.
$$
If $\mathcal{M}$ is a pseudo-Riemannian manifold, set
$$\mathcal{U}_{\mathcal{M},P}^k:=(T_PM,g_{\mathcal{M}}|_{T_PM},R_{\mathcal{M}}|_{T_PM},...,
      \nabla^kR_{\mathcal{M}}|_{T_PM})\,.
$$
We define $\mathcal{U}^\infty$ and $\mathcal{U}_{\mathcal{M},P}^\infty$ similarly using an infinite sequence.
One says that $\mathcal{U}^k$ is a {\it $k$-model for
$\mathcal{M}$} if for every point $P$ of $M$, there is a isomorphism $\Phi_P$ from
$\mathcal{U}_{\mathcal{M},P}^k$ to $\mathcal{U}^k$, i.e. $\Phi_P:T_PM\rightarrow V$ has
$$\Phi_P^*\langle\cdot,\cdot\rangle=g_{\mathcal{M}}|_{T_PM}\quad\text{and}\quad
  \Phi_P^*A^i=\nabla^iR_{\mathcal{M}}|_{T_PM}\quad\text{for}\quad0\le i\le k\,.
$$
One says that $\mathcal{M}$ is {\it $k$-curvature homogeneous} if it admits a $k$-model; this means that the metric and the
covariant derivatives of the curvature up to order $k$ ``look the same at every point''. One says that
$\mathcal{M}$ is
{\it $k$-modeled on a homogeneous space $\mathcal{N}$} if there exists $Q\in N$ so that $\mathcal{U}_{\mathcal{N},Q}^k$
is a $k$-model for $\mathcal{M}$. Let
\begin{eqnarray*}
&&\mathcal{C}_2(\mathcal{O}):=\{f\in\mathcal{A}(\mathcal{O}):f^{(2)}>0\},\\
&&\mathcal{C}_3(\mathcal{O}):=\{f\in\mathcal{A}(\mathcal{O}):f^{(2)}>0\text{ and }f^{(3)}>0\},\\
&&\alpha_p(f):=f^{(p+2)}\{f^{(2)}\}^{p-1}\{f^{(3)}\}^{-p}\quad\text{for }f\in\mathcal{C}_3(\mathcal{O})\text{ and }p\ge2\,.
\end{eqnarray*}

We extend an earlier result of Derdzinski:

\begin{theorem}\label{thm-1.6}
\ \begin{enumerate}
\item If $f\in\mathcal{C}_2(\mathbb{R})$, $\mathcal{M}_f$ is $0$-modeled on the irreducible symmetric space $\mathcal{M}_{y^2}$.
\item If $f\in\mathcal{C}_3(\mathbb{R})$, $\mathcal{M}_f$ is $1$-modeled on the homogeneous space $\mathcal{M}_{e^y}$.
\item Let $f_i\in\mathcal{C}_3(\mathbb{R})$. The following assertions are equivalent:
\begin{enumerate}\item There exists an isometry
$\phi:(\mathcal{M}_{f_1},P_1)\rightarrow(\mathcal{M}_{f_2},P_2)$.
\item We have
$\alpha_p(f_1)(P_1)=\alpha_p(f_2)(P_2)$ for $p\ge2$.
\end{enumerate}\end{enumerate}\end{theorem}

\begin{remark}\label{thm-1.7}\rm The manifold $\mathcal{M}_{e^y+e^{2y}}$ is a complete pseudo-Riemannian manifold of signature
$(2,2)$ which is not homogeneous but which is $0$-modeled on the irreducible symmetric space $\mathcal{M}_{y^2}$. Work of
Tricerri and Vanhecke
\cite{TV86} shows this is not possible in the Riemannian setting; work of Cahen, Leroy, Parker, Tricerri, and Vanhecke \cite{CLPT90}
shows this is not possible in the Lorentzian setting. 
\end{remark}

Although the focus of this paper is primarily on global questions, we  can also discuss the local geometry. If
$f\in\mathcal{A}(\mathcal{O})$, set
$$\mathcal{M}_{f,\mathcal{O}}:=(\{(x,y,\tilde x,\tilde y)\in\mathbb{R}^4:y\in\mathcal{O}\},g_f)\,.$$
\begin{theorem}\label{thm-1.8}
 Let $f\in\mathcal{C}_3(\mathcal{O})$.
The following assertions are equivalent:
\begin{enumerate}\item $\mathcal{M}_{f,\mathcal{O}}$ is locally homogeneous.
\item $\mathcal{M}_{f,\mathcal{O}}$ is $2$-curvature homogeneous.
\item $\alpha_2$ is constant on $\mathcal{O}$.
\item either $f^{(2)}=a(y+b)^c$ or $f=ae^\lambda y$.
\end{enumerate}\end{theorem}

\begin{remark}\label{rmk-1.9}\rm The manifold $\mathcal{M}_{y^4}$ is a
connected complete pseudo-Riemannian manifold of signature $(2,2)$ which contains a proper open connected homogeneous submanifold
$\mathcal{M}_{y^4,\mathbb{R}^+}$. This is not possible in the Riemannian
setting. 
\end{remark}
\begin{remark}\label{rmk-1.10}\rm Work of Singer \cite{S60} in the Riemannian setting and of Podesta and Spiro \cite{PS04} in the
higher signature setting shows there exists a universal integer $k_{p,q}$ so that if $\mathcal{M}$ is a complete simply connected
$k_{p,q}$-curvature homogeneous manifold of signature
$(p,q)$, then $\mathcal{M}$ is in fact homogeneous. Opozda
\cite{O97} has established a similar result in the affine setting. Gromov \cite{Gr86} and Yamato \cite{Y89} have given upper
bounds for $k_{0,q}$ which are linear in $q$. There are Riemannian manifolds which are
$0$-curvature homogeneous but not homogeneous, see for example Ferus, Karcher, and M\"unzer
\cite{FKM81} or Takagi \cite{T74}. It is clear that $k_{0,2}=k_{1,1}=k_{2,0}=0$. Work of Sekigawa, Suga, and Vanhecke
\cite{SSV92,SSV95} shows
$k_{0,3}=k_{0,4}=1$. We refer to the discussion in Boeckx, Kowalski, and Vanhecke
\cite{BKV96} for further details concerning $k$-curvature homogeneous manifolds in the Riemannian setting. In the Lorentzian
setting, work of Bueken and Djori\'c \cite{BD00} and of Bueken and Vanhecke \cite{BV97} shows that $k_{1,2}\ge2$.
Theorem \ref{thm-1.8} shows that $k_{2,2}\ge2$; furthermore, $2$-curvature homogeneity implies local homogeneity for this family. For
$p\ge3$, there are manifolds of neutral signature $(p,p)$ which are $(p-1)$-curvature homogeneous but which are not
homogeneous \cite{GN04b}; these examples have much the same flavor as the manifolds $\mathcal{M}_f$ discussed here.\end{remark}

We now return to the global setting. Let $G_f$ be the Lie group of isometries of $\mathcal{M}_f$ and let $\gg_f$ be the associated
Lie algebra. As we are working in the analytic category, Theorem \ref{thm-1.3} (3) shows that any local isometry extends to a
global isometry so we may identify $\gg_f$ with the Lie algebra of Killing vector fields on
$\mathcal{M}_f$. 
\begin{theorem}\label{thm-1.11}\ \begin{enumerate}
\item If $f^{(2)}=0$, then $\dim\gg_f=10$.
\item If $f^{(2)}=c\ne0$, then $\dim\gg_f=8$.
\item If $f^{(2)}=ae^{\lambda y}$ for $a\ne0$ and $\lambda\ne0$, then $\dim\gg_f=6$.
\item If $f^{(2)}=a(y+b)^n$ for $a\ne0$ and $n=1,2,...$, then $\dim\gg_f=6$.
\item If $f^{(2)}\ne ae^{\lambda y}$ and if $f^{(2)}\ne a(y+b)^n$ for $n\in\mathbb{N}$, then $\dim\gg_f=5$.
\end{enumerate}
\end{theorem}

\begin{remark}\label{rmk-1.12}\rm 
The structure of the Killing vector fields on a manifold reflects the underlying geometry of the manifold. The manifolds
described in Theorem
\ref{thm-1.11} (1) are flat, those in (2) are symmetric but not flat, those in (3)
are homogeneous, those in (4) are homogeneous on the open sets where $y+b$ is positive or negative but are not globally homogeneous,
and those in (5) are not homogeneous. The calculations we shall perform in Section \ref{sect-6} not only determines the dimensions,
but also exhibits bases explicitly for these algebras and permit a determination of their structure constants.\end{remark}

Here is a brief outline to this paper. In Section \ref{sect-2}, we determine the geodesics and the curvature tensor of the manifolds
$\mathcal{M}_f$. Theorem \ref{thm-1.3} (1)-(4) and Theorem \ref{thm-1.4} (1) follow. In Section \ref{sect-3}, we study models and
establish Assertions (5) and (6) of Theorem
\ref{thm-1.3} and Assertions (1) and (2) of Theorem \ref{thm-1.6}. Section \ref{sect-4} deals with isometries. In the real
analytic context, an isomorphism from $\mathcal{U}_{\mathcal{M},P}^\infty$ to $\mathcal{U}_{\mathcal{N},Q}^\infty$
induces a local isometry from $\mathcal{M}$ to $\mathcal{N}$. Since the manifolds we are considering are simply connected and
complete, the local isometry extends globally. We use this observation to prove Theorem \ref{thm-1.3} (7) and Theorem \ref{thm-1.6}
(3). We show that the solutions to the differential equation $hh^{\prime\prime}=kh^\prime h^\prime$ have the form $h=ae^{\lambda y}$
if
$k=1$ and $h=a(y+b)^c$ for $c=(1-k)^{-1}$ if $k\ne1$. This observation is used to establish Theorem
\ref{thm-1.4} (2) and Theorem \ref{thm-1.8}.
In the final two sections of the paper, we study the Killing vector fields. In Section
\ref{sect-5}, we establish a structure theorem for Killing vector fields on generalized plane wave manifolds satisfying a
non-degeneracy condition; we then specialize this result to the manifolds $\mathcal{M}_f$ when $f^{(2)}\ne0$. In Section
\ref{sect-6}, we prove Theorem \ref{thm-1.11}.

Throughout this paper, we shall be studying the case when $f^{(2)}>0$; the case $f^{(2)}<0$ is entirely analogous. The sign of
$f^{(3)}$ is irrelevant; the crucial condition is that $f^{(3)}$ is non-zero.

\section{Geodesics and Curvature}\label{sect-2}
\begin{definition}\label{defn-2.1}\rm We follow the discussion in \cite{GN04a}. For $p\ge2$, let $(x_1,...,x_p,\tilde x_1,...,\tilde
x_p)$ be coordinates on
$\mathbb{R}^{2p}$. Let indices $i$, $j$, $k$ range from $1$ through $p$. Set 
$$\partial_i^x:=\ffrac{\partial}{\partial x_i}\quad\text{and}\quad
  \partial_i^{\tilde x}:=\ffrac{\partial}{\partial\tilde x_i}\,.
$$
Let $\Xi:=\Xi_{ij}(x_1,...,x_p)$ be a smooth symmetric $2$-tensor field on $\mathbb{R}^p$. We consider the generalized plane
wave manifold
$\mathcal{P}_\Xi:=(\mathbb{R}^{2p},g_\Xi)$ of signature $(p,p)$ where
$$
g_\Xi(\partial_i^x,\partial_j^x):=\Xi_{ij}(x_1,...,x_p)\quad\text{and}\quad 
g_\Xi(\partial_i^x,\partial_i^{\tilde x}):=\delta_{ij}\,.
$$\end{definition}

\begin{lemma}\label{lem-2.2}Let $\Gamma_{ijk}^x:=\ffrac12(\partial_i^x\Xi_{jk}+\partial_j^x\Xi_{ik}-\partial_k^x\Xi_{ij})$. 
\begin{enumerate}
\item The non-zero components
of $\nabla^kR_{\mathcal{P}_\Xi}$ are
$$\nabla^kR_{\mathcal{P}_\Xi}(\partial_{i_1}^x,\partial_{i_2}^x,\partial_{i_3}^x,\partial_{i_4}^x;\partial_{j_1}^x,...,\partial_{j_k}^x)
    =(\partial_{j_1}^x...\partial_{j_k}^x)(\partial_{i_1}^x\Gamma_{i_2i_3i_4}^x-\partial_{i_2}^x\Gamma_{i_1i_3i_4}^x)\,.$$
\item All geodesics in $\mathcal{P}_\Xi$ extend for infinite time.
\item If $P\in\mathbb{R}^{2p}$, then $\exp_{\mathcal{P}_\Xi,P}:T_P\mathbb{R}^4\rightarrow\mathbb{R}^4$ is a diffeomorphism.
\item If $\Xi$ is quadratic in $(x_1,...,x_p)$, then $\mathcal{M}_{\Xi}$ is a symmetric space.
\item All the local scalar Weyl invariants of $\mathcal{P}_\Xi$ vanish.
\end{enumerate}\end{lemma}
\begin{proof} The non-zero Christoffel symbols are:
\begin{equation}\label{eqn-2.a}
g_\Xi(\nabla_{\partial_i^x}\partial_j^x,\partial_k^x)=\Gamma_{ijk}^x\quad\text{and}\quad
\nabla_{\partial_i^x}\partial_j^x=\textstyle\sum_k\Gamma_{ijk}^x\partial_k^{\tilde x}\,.
\end{equation}
Assertion (1) follows by direct computation. We use Equation (\ref{eqn-2.a}) to see that the geodesics have the form
\begin{eqnarray*}
&&x_k(t)=\alpha_k+\beta_kt,\\
&&\tilde x_k(t)=\tilde\alpha_k+\tilde\beta_kt
-\textstyle\sum_{ij}\beta_i\beta_j\int_{s=0}^t\int_{r=0}^s\Gamma_{ijk}(x(r))drds,
\end{eqnarray*}
they extend for all time. Furthermore, given $P$ and $Q$ in $\mathbb{R}^{2p}$, there is a
unique geodesic $\sigma$ with $\sigma(0)=P$ and $\sigma(1)=Q$; thus $\exp_{\mathcal{P}_\Xi,P}$ is a diffeomorphism from
$T_P\mathbb{R}^{2p}$ to
$\mathbb{R}^{2p}$; we refer to \cite{GN04a} for further details. This establishes Assertions (2) and (3); Assertion (4) is an
immediate consequence of Assertion (1).

Introduce a new frame 
\begin{equation}\label{eqn-2.d}
X_i:=\partial_i^x-\ffrac12\textstyle\sum_j\Xi_{ij}\partial_j^{\tilde x}\quad\text{and}\quad
\tilde X_i:=\partial_i^{\tilde x}\,.
\end{equation}
Then $\{X_i,\tilde X_i\}$ is a hyperbolic frame, i.e. the only non-zero components of the metric tensor are given by
$g_\Xi(X_i,\tilde X_j)=\delta_{ij}$. Assertion (5) follows since $\nabla^kR_{\mathcal{P}_\Xi}$ is supported on the
totally isotropic space $\operatorname{Span}\{X_i\}$ for any $k$.
\end{proof}

\begin{proof}[Proof of Theorem \ref{thm-1.3} (1)-(4) and Theorem \ref{thm-1.4} (1)] We have:
$$\begin{array}{ll}
g(\nabla_{\partial x}\partial_x,\partial_y)=\partial_yf,&
g(\nabla_{\partial x}\partial_y,\partial_x)=g(\nabla_{\partial_y}\partial_x,\partial_x)=-\partial_yf,\\
\nabla_{\partial_x}\partial_x=(\partial_yf)\partial_{\tilde y},&
\nabla_{\partial_x}\partial_y=\nabla_{\partial_y}\partial_x=-(\partial_yf)\partial_{\tilde x}.
\vphantom{\vrule  height 12pt}
\end{array}$$
Assertion (1) of Theorem \ref{thm-1.3} now follows. Assertions (2)-(4) of Theorem \ref{thm-1.3} and Assertion (1) of
Theorem \ref{thm-1.4} follow directly from Lemma \ref{lem-2.2}.\end{proof}

We shall follow the discussion in \cite{DG04} to discuss the following hypersurfaces.
\begin{definition}\label{defn-2.3}\rm
Give
$\mathbb{R}^{p,p+1}:=\operatorname{Span}\{e_1,...,e_p,\tilde e_1,...,\tilde e_p,\check e\}$ the inner-product:
$$\langle e_i,\tilde e_j\rangle=\delta_{ij},\quad\langle \check e,\check e\rangle=1\,.$$
Let $\psi$ be a smooth function on $\mathbb{R}^p$. Define an embedding $\Psi:\mathbb{R}^{2p}\rightarrow\mathbb{R}^{(p,p+1)}$ by:
$$\Psi(x_1,...,x_p,\tilde x_1,...,\tilde x_p):=\sum_{i=1}^p\bigg\{x_ie_i+\tilde x_i\tilde e_i\bigg\}+\psi(x_1,...,x_p)\check e\,.$$
Let
$h_\psi:=\Psi^*\langle\cdot,\cdot\rangle$ and let
$\mathcal{H}_\psi:=(\mathbb{R}^4,h_\psi)$ be the associated pseudo-Riemannian manifold of signature $(2,2)$.
\end{definition}

\begin{lemma}\label{lem-2.4}Let $L_{ij}:=\partial_i^x\partial_j^x\psi$. \begin{enumerate}
\item  The non-zero components of $\nabla^kR_{\mathcal{H}_\psi}$ are:
$$\nabla^kR_{\mathcal{H}_\psi}(\partial_{i_1}^x,\partial_{i_2}^x,\partial_{i_3}^x,\partial_{i_4}^x;\partial_{j_1}^x,...,\partial_{j_k}^x)
    =(\partial_{j_1}^x...\partial_{j_k}^x)(L_{i_1i_4}L_{i_2i_3}-L_{i_1i_3}L_{i_2i_4})\,.$$
\item All geodesics in $\mathcal{H}_\psi$ extend for infinite time.
\item If $P\in\mathbb{R}^{2p}$, then $\exp_{\mathcal{H}_\psi,P}:T_P\mathbb{R}^4\rightarrow\mathbb{R}^4$ is a diffeomorphism.
\item If $\psi(x_1,x_2)=\frac12x_1^2+f(x_2)$, then the non-zero components of $\nabla^kR_{\mathcal{H}_\psi}$ are
$$\nabla^kR_{\mathcal{H}_\psi}(\partial_1^x,\partial_2^x,\partial_2^x,\partial_1^x;\partial_2^x,...,\partial_2^x)=f^{(2+k)}\,.$$
\end{enumerate}
\end{lemma}

\begin{proof}
Since $\mathcal{H}_\psi=\mathcal{P}_\Xi$ where $\Xi_{ij}=\partial_i^x\psi\cdot\partial_j^x\psi$, $\mathcal{H}_\psi$ is a manifold of
the form discussed above. Consequently, Assertion (1) can be derived from Lemma \ref{lem-2.2}. It is, however, instructive to compute
this directly. As
$$\nu:=-\partial_1^x\psi\tilde e_1-...-\partial_p^x\psi\tilde e_p+\check e$$
is the unit normal to the hypersurface,
$L_{ij}$ is the second fundamental form of the embedding. This establishes Assertion (1) for $k=0$. The computation
of
$\nabla^kR_{\mathcal{H}_\psi}$ for
$k>0$ now follows from the structure equations given in Equation (\ref{eqn-2.a}); the Christoffel symbols play no role.
Assertions (2) and (3) follow from Lemma \ref{lem-2.2}. Assertion (4) follows from Assertion (1).\end{proof}

\section{Models}\label{sect-3}
\begin{definition}\label{defn-3.1}\rm Let $\{X,Y,\tilde X,\tilde Y\}$ be a basis for $\mathbb{R}^4$. Let
$\mathcal{U}^0:=(\mathbb{R}^4,\langle\cdot,\cdot\rangle,A^0)$ and $\mathcal{U}^1:=(\mathbb{R}^4,\langle\cdot,\cdot\rangle,A^0,A^1)$
where the non-zero components of
$\langle\cdot,\cdot\rangle$, $A^0$, and $A^1$ are:
$$\langle X,\tilde X\rangle=\langle Y,\tilde Y\rangle=1,\ \ A^0(X,Y,Y,X)=1,\ \ \text{and}\ \ A^1(X,Y,Y,X;Y)=1\,.$$
\end{definition}
\begin{proof}[Proof of Theorem \ref{thm-1.6} (1)] 
Set $\partial_x:=\frac{\partial}{\partial x}$, $\partial_y:=\frac{\partial}{\partial y}$, $\partial_{\tilde
x}:=\frac{\partial}{\partial \tilde x}$, and
$\partial_{\tilde y}:=\frac{\partial}{\partial\tilde y}$. Then:
$$g_f(\partial_x,\partial_x)=-2f,\quad g_f(\partial_x,\partial_{\tilde x})=g_f(\partial_y,\partial_{\tilde y})=1,\quad
R_{\mathcal{M}_f}(\partial_x,\partial_y,\partial_y,\partial_x)=f^{(2)}\,.$$
Suppose $f^{(2)}>0$. To see that $\mathcal{U}^0$ is $0$-curvature model for $\mathcal{M}_f$, we set
$$
X:=\partial_x+f\partial_{\tilde x},\quad
  Y:=(f^{(2)})^{-1/2}\partial_y,\quad \tilde X:=\partial_{\tilde x},\quad\tilde Y:=\partial_{\tilde y}\,.
$$
Assertion (1) of Theorem \ref{thm-1.6} now follows.
\end{proof}

\begin{proof}[Proof of Theorem \ref{thm-1.3} (5, 6)] In view of the above discussion, to prove Theorem \ref{thm-1.3} (5), it suffices
to show that
$\mathcal{U}^0$ is spacelike and timelike Jordan Osserman. Let
$\xi=aX+bY+\tilde a\tilde X+\tilde b\tilde Y\in\mathbb{R}^4$. If
$\xi$ is not null, then
$(a,b)\ne(0,0)$. We compute:
$$\begin{array}{ll}
R_{\mathcal{U}^0}(X,Y)X=-\tilde Y,&R_{\mathcal{U}^0}(X,Y)Y=\tilde X,\\
J_{\mathcal{U}^0}(\xi)(aX+bY)=0,&J_{\mathcal{U}^0}(\xi)\tilde X=0,\vphantom{\vrule height 12pt}\\ 
J_{\mathcal{U}^0}(\xi)(-bX+aY)=(a^2+b^2)(-b\tilde X+a\tilde Y),&J_{\mathcal{U}^0}(\xi)\tilde Y=0\,.\vphantom{\vrule height 12pt}
\end{array}$$ 
Thus $J_{\mathcal{U}^0}(\xi)^2=0$ and $\operatorname{rank}\{J_{\mathcal{U}^0}(\xi)\}=1$. This shows that the Jordan normal form of
$J_{\mathcal{U}^0}(\cdot)$ is constant on the set of non-null vectors and hence $\mathcal{U}^0$ is spacelike and timelike Jordan
Osserman.

Similarly, to establish Theorem \ref{thm-1.3} (6), it suffices to show that
$\mathcal{U}^0$ is spacelike and timelike Ivanov-Petrova.
Let $\{e_1,e_2\}$ be an oriented orthonormal basis for an oriented $2$ plane $\pi$ which contains no non-zero null vectors. We expand
$$e_1=a_1X+b_1Y+\tilde a_1\tilde X+\tilde b_1\tilde Y\quad\text{and}\quad e_2=a_2X+b_2Y+\tilde a_2\tilde X+\tilde b_2\tilde Y$$
where $a_1b_2-a_2b_1\ne0$. We have
$R_{\mathcal{U}^0}(\pi)=(a_1b_2-a_2b_1)R_{\mathcal{U}^0}(X,Y)$. Thus
$$\begin{array}{ll}
R_{\mathcal{U}^0}(\pi):X\rightarrow-(a_1b_2-a_2b_1)\tilde Y,&
  R_{\mathcal{U}^0}(\pi):Y\rightarrow(a_1b_2-a_2b_1)\tilde X,\\
  R_{\mathcal{U}^0}(\pi):\tilde X\rightarrow0,&
  R_{\mathcal{U}^0}(\pi):\tilde Y\rightarrow0\,.
\end{array}$$
Consequently,
$R_{\mathcal{U}^0}(\pi)^2=0$ and $\operatorname{rank}\{R_{\mathcal{U}^0}(\pi)\}=2$ so $\mathcal{U}^0$ is spacelike and timelike
Jordan Ivanov-Petrova.
\end{proof}

Let $\mathcal{U}_{f,P}^\infty:=(\mathbb{R}^4,\langle\cdot,\cdot\rangle,A_{f,P}^0,...)$
where
\begin{equation}\label{eqn-3.a}
\langle X,\tilde X\rangle=\langle Y,\tilde Y\rangle:=1,\quad\text{and}\quad
A_{f,P}^k(X,Y,Y,X;Y,...,Y):=f^{(k+2)}(P)\,.
\end{equation}
If $f\in\mathcal{C}_3$, let $\mathcal{V}_{f,P}^\infty:=(\mathbb{R}^4,\langle\cdot,\cdot\rangle,B_{f,P}^0,...)$ where
\begin{eqnarray*}
&&\langle X,\tilde X\rangle=\langle Y,\tilde Y\rangle:=1,\quad B_{f,P}^0(X,Y,Y,X)=1,\\
&&B_{f,P}^1(X,Y,Y,X;Y)=1,\quad\text{and}\\
&&B_{f,P}^k(X,Y,Y,X;Y,...,Y):=\alpha_k(f,P)\quad\text{for}\quad k\ge2\,.
\end{eqnarray*}

\begin{lemma}\label{lem-3.2}\ \begin{enumerate}
\item There exists an isomorphism between $\mathcal{U}_{f,P}^\infty$ and $\mathcal{U}_{\mathcal{M}_f,P}^\infty$.
\item There exists an isomorphism between $\mathcal{U}_{f,P}^\infty$ and $\mathcal{U}_{\mathcal{H}_\psi,P}^\infty$ for
$\psi=\frac12x_1^2+f(x_2)$.
\item If $f\in\mathcal{C}_3(\mathbb{R})$, then there exists an isomorphism between $\mathcal{U}_{f,P}^\infty$ and
$\mathcal{V}_{f,P}^\infty$.
\end{enumerate}
\end{lemma}

\begin{proof} To prove Assertion (1), we set
$$X:=\partial_x+f\partial_{\tilde x},\quad
  Y:=\partial_y,\quad \tilde X:=\partial_{\tilde x},\quad\tilde Y:=\partial_{\tilde y}\,.
$$
This normalizes the metric to be hyperbolic but does not change the curvature tensor; the relations of Equation (\ref{eqn-3.a}) then
hold by Theorem \ref{thm-1.3} (1). To prove Assertion (2), we set
$$\begin{array}{ll}
X:=\partial_1^x-\frac12h_\psi(\partial_1^x,\partial_1^x)\partial_1^{\tilde x}
  -\frac12h_\psi(\partial_1^x,\partial_2^x)\partial_2^{\tilde x},&
\tilde X:=\partial_1^{\tilde x},\\
Y:=\partial_2^x-\frac12h_\psi(\partial_2^x,\partial_1^x)\partial_1^{\tilde x}
  -\frac12h_\psi(\partial_2^x,\partial_2^x)\partial_2^{\tilde x},&
\tilde Y:=\partial_2^{\tilde x}\,.\vphantom{\vrule height 12pt}
\end{array}$$
Again, this normalizes the metric to be  hyperbolic but does not change the curvature tensor; the relations of Equation
(\ref{eqn-3.a}) then hold by Lemma \ref{lem-2.4}.

Let $f\in\mathcal{C}_3(\mathbb{R})$. Set $X_1:=\varepsilon_1X$, $Y_1:=\varepsilon_2Y$, $\tilde X_1:=\varepsilon_1^{-1}\tilde X$,
$\tilde Y_1=\varepsilon_2^{-1}\tilde Y$ to define a hyperbolic basis with
$$\nabla^kR(X_1,Y_1,Y_1,X_1;Y_1,...,Y_1)=\varepsilon_1^2\varepsilon_2^{k+2}f^{(k+2)}\,.$$
To ensure that
$$A_{f,P}^0(X_1,Y_1,Y_1,X_1)=1\quad\text{and}\quad A_{f,P}^1(X_1,Y_1,Y_1,X_1;Y_1)=1,$$
we must have that
$\varepsilon_1^2\varepsilon_2^2f^{(2)}=1$ and that $\varepsilon_1^2\varepsilon_2^3f^{(3)}=1$. We set
$$\varepsilon_2:=\ffrac{f^{(2)}}{f^{(3)}}\quad\text{and}\quad\varepsilon_1:=\varepsilon_2^{-1}\{f^{(2)}\}^{-1/2}\,.$$
This then yields for $k\ge2$ that
\begin{eqnarray*}
&&\nabla^kR(X_1,Y_1,Y_1,X_1;Y_1,...,Y_1)=\varepsilon_2^{-2}\{f^{(2)}\}^{-1}\varepsilon_2^{2+k}f^{(2+k)}\\
&=&\{f^{(2)}\}^{k-1}\{f^{(3)}\}^{-k}f^{(2+k)}=\alpha_k(f)
\,.
\end{eqnarray*}
This establishes the desired isomorphism.
\end{proof}

\begin{proof}[Proof of Theorem \ref{thm-1.6} (2)] If $f\in\mathcal{C}_3(\mathbb{R})$, then we may restrict the isomorphism of Lemma
\ref{lem-3.2} to see that there is an isomorphism between
$\mathcal{U}_{\mathcal{M}_f,P}^1$ and  $\mathcal{U}_{f,P}^1$ and that there is an isomorphism between $\mathcal{U}_{f,P}^1$ and
$\mathcal{V}_{f,P}^1=\mathcal{U}^1$. Because $\mathcal{U}^1$ is depends neither on $P$ nor on $f$, this shows that $\mathcal{M}_f$ is
$1$-curvature homogeneous and is, in particular, modeled on $\mathcal{M}_{e^y}$.
\end{proof}

We shall need the following technical Lemma in Section \ref{sect-6}:

\begin{lemma}\label{lem-3.3}\ \begin{enumerate}\item
If $G_0$ is the symmetry group of $\mathcal{U}^0$, then $G_0\subset GL(4,\mathbb{R})$ is the $4$ dimensional
Lie group with 2 connected components described by:
$$G_0=\left\{\left(\begin{array}{ll}\alpha&\gamma(\alpha^{-1})^t\\0&(\alpha^{-1})^t\end{array}\right):\alpha,\gamma\in
M_2(\mathbb{R}),\ \det\alpha=\pm1\quad\text{and}\quad\gamma+\gamma^t=0\right\}\,.$$
\item If $G_1$ is the symmetry group of $\mathcal{U}^1$, then $G_1\subset GL(4,\mathbb{R})$ is the $2$ dimensional connected
Lie group described by:
$$G_1=\left\{\left(\begin{array}{ll}\alpha&\gamma(\alpha^{-1})^t\\0&(\alpha^{-1})^t\end{array}\right):\alpha,\gamma\in
M_2(\mathbb{R}),\ \alpha=
\left(\begin{array}{ll}1&0\\a_{21}&1\end{array}\right)\text{ and }\gamma+\gamma^t=0\right\}\,.$$

\end{enumerate}\end{lemma}

\begin{proof}
If $\Theta\in G_0$, let
\begin{eqnarray*}
&&\Theta X=a_{11}X+a_{12}Y+a_{13}\tilde X+a_{14}\tilde Y,\\
&&\Theta Y=a_{21}X+a_{22}Y+a_{23}\tilde X+a_{24}\tilde Y,\\
&&\Theta\tilde X=a_{31}X+a_{32}Y+a_{33}\tilde X+a_{34}\tilde Y,\\
&&\Theta\tilde Y=a_{41}X+a_{42}Y+a_{43}\tilde X+a_{44}\tilde Y\,.
\end{eqnarray*}
We set
\begin{eqnarray*}
&&\Theta=\left(\begin{array}{llll}
a_{11}&a_{12}&a_{13}&a_{14}\\a_{21}&a_{22}&a_{23}&a_{24}\\
a_{31}&a_{32}&a_{33}&a_{34}\\a_{41}&a_{42}&a_{43}&a_{44}\end{array}\right)=
\left(\begin{array}{ll}\alpha_1&\alpha_2\\\alpha_3&\alpha_4\end{array}\right)\quad\text{where}\\
&&\alpha_1:=\left(\begin{array}{ll}a_{11}&a_{12}\\a_{21}&a_{22}\end{array}\right),\quad
\alpha_2:=\left(\begin{array}{ll}a_{13}&a_{14}\\a_{23}&a_{24}\end{array}\right)\\
&&\alpha_3:=\left(\begin{array}{ll}a_{31}&a_{32}\\a_{41}&a_{42}\end{array}\right),\quad
\alpha_4:=\left(\begin{array}{ll}a_{33}&a_{34}\\a_{43}&a_{44}\end{array}\right)\,.
\end{eqnarray*}
As $R(\Theta X,\Theta Y,\Theta Y,\Theta X)=1$, we have $(a_{11}a_{22}-a_{12}a_{21})^2=1$ and $\det(\alpha_1)^2=1$. As
\begin{eqnarray*}
0&=&R(\Theta X,\Theta Y,\Theta X,\Theta \tilde X)=R(\Theta X,\Theta Y,\Theta Y,\Theta\tilde X)\\
&=&R(\Theta X,\Theta Y,\Theta X,\Theta \tilde Y)=R(\Theta X,\Theta Y,\Theta Y,\Theta\tilde Y),
\end{eqnarray*}
we have that $a_{31}=a_{32}=a_{41}=a_{42}=0$ and thus $\alpha_3=0$.
As 
$$\begin{array}{ll}
g(\Theta X,\Theta\tilde X)=1,&g(\Theta X,\Theta\tilde Y)=0\\
g(\Theta Y,\Theta\tilde X)=0,&g(\Theta Y,\Theta\tilde Y)=1,\vphantom{\vrule height 11pt}\end{array}
$$
we have the relations
$$\begin{array}{ll}
1=a_{11}a_{33}+a_{12}a_{34},&0=a_{11}a_{43}+a_{12}a_{44},\\
0=a_{21}a_{33}+a_{22}a_{34},&1=a_{21}a_{43}+a_{22}a_{44}\vphantom{\vrule height 12pt}
\end{array}$$
which can be rewritten in matrix form:
$$\id=\left(\begin{array}{ll}a_{11}&a_{12}\\a_{21}&a_{22}\end{array}\right)
    \left(\begin{array}{ll}a_{33}&a_{43}\\a_{34}&a_{44}\end{array}\right)\quad\text{i.e.}\quad
   I=\alpha_1\alpha_4^t\,.$$
The relations $g(\Theta X,\Theta X)=g(\Theta X,\Theta Y)=g(\Theta Y,\Theta Y)=0$ yield the equations
$$\begin{array}{l}
0=a_{11}a_{13}+a_{12}a_{14},\\
0=a_{11}a_{23}+a_{12}a_{24}+a_{13}a_{21}+a_{14}a_{22},\\
0=a_{21}a_{23}+a_{22}a_{24}
\end{array}$$
which can be rewritten in matrix form as
$$0=\left(\begin{array}{ll}a_{11}&a_{12}\\a_{21}&a_{22}\end{array}\right)
\left(\begin{array}{ll}a_{13}&a_{23}\\a_{14}&a_{24}\end{array}\right)
+\left(\begin{array}{ll}a_{13}&a_{14}\\a_{23}&a_{24}\end{array}\right)
\left(\begin{array}{ll}a_{11}&a_{21}\\a_{12}&a_{22}\end{array}\right)$$
or equivalently as
$\alpha_1\alpha_2^t+\alpha_2\alpha_1^t=0$.  Setting $\alpha=\alpha_1$ and $\alpha_2=\gamma(\alpha^{-1})^t$ yields
$\gamma+\gamma^t=0$ which establishes the first implication of Assertion (1). Conversely, the implications are all reversible and
thus $\Theta$ satisfies these relations implies $\Theta\in G_0$.

Suppose $\Theta\in G_1\subset G_0$. Since 
\begin{eqnarray*}
&&R(\Theta X,\Theta Y,\Theta Y,\Theta X)=1,\\ 
&&\nabla R(\Theta X,\Theta Y,\Theta Y,\Theta X;\Theta Y)=1,\\
&&\nabla R(\Theta X,\Theta Y,\Theta Y,\Theta X;\Theta X)=0,
\end{eqnarray*}
we have $a_{11}a_{22}-a_{12}a_{21}=1$, $a_{22}=1$, and $a_{12}=0$.
Thus we have $a_{11}=1$ as well and $\alpha$ has the desired form.
\end{proof}

\section{Isometries}\label{sect-4}
We shall need the following useful observation:
\begin{lemma}\label{lem-4.1} Let $\mathcal{M}_i:=(M_i,g_i)$ be real analytic pseudo-Riemannian manifolds for $i=1,2$. Assume
there exist points $P_i\in M_i$ so
$\exp_{\mathcal{M}_i,P_i}:T_{P_i}M_i\rightarrow M_i$ is a diffeomorphism and so there exists an isomorphism $\Phi$ between
$\mathcal{U}_{\mathcal{M}_1,P_1}^\infty$ and $\mathcal{U}_{\mathcal{M}_2,P_2}^\infty$. Then
$\phi:=\exp_{\mathcal{M}_2,P_2}\circ\Phi\circ\exp_{\mathcal{M}_1,P_1}^{-1}$ is an isometry from $\mathcal{M}_1$ to $\mathcal{M}_2$.
\end{lemma}

\begin{proof} Belger and Kowalski \cite{BeKo94} note about analytic pseudo-Riemannian metrics that the
``metric
$g$ is uniquely determined, up to local isometry, by the tensors $R$, $\nabla R$, ..., $\nabla^kR$, ... at one point.''; see also
Gray \cite{Gr73} for related work. The desired result now follows.
\end{proof}

\begin{proof}[Proof of Theorem \ref{thm-1.3} (7)] Suppose $f\in\mathcal{C}_3(\mathbb{R})$. Let $\psi(x_1,x_2)=\frac12x_1^2+f(x_2)$.
By Lemma
\ref{lem-3.2},
$\mathcal{U}_{\mathcal{M}_f,P}^\infty$ and
$\mathcal{U}_{\mathcal{H}_\psi,P}^\infty$ are isomorphic. By Theorem \ref{thm-1.3} (3) and Lemma \ref{lem-2.4} (3), the exponential
map is a global diffeomorphism for both manifolds. The desired result now follows by Lemma \ref{lem-4.1}.\end{proof}

\begin{proof}[Proof of Theorem \ref{thm-1.6} (3)] Suppose $f^{(3)}(P)\ne0$. We may then choose $X$ and $Y$ in $T_P\mathbb{R}^4$ so
that $\nabla R_{\mathcal{M}_f}(X,Y,Y,X;Y)\ne0$; for example, we could set $X=\partial_x$ and $Y=\partial_y$. We expand
$$X=a_1\partial_x+a_2\partial_y+\tilde a_1\partial_{\tilde x}+\tilde a_2\partial_{\tilde y}\quad\text{and}\quad
Y=b_1\partial_x+b_2\partial_y+\tilde b_1\partial_{\tilde x}+\tilde b_2\partial_{\tilde y}\,.
$$
As
$\nabla R_{\mathcal{M}_f}(X,Y,Y,X;Y)=(a_1b_2-a_2b_1)^2b_2f^{(3)}$, $(a_1b_2-a_2b_1)^2b_2\ne0$ and
\begin{eqnarray*}
&&\nabla^kR_{\mathcal{M}_f}(X,Y,Y,X;Y,...,Y)=(a_1b_2-a_2b_1)^2b_2^kf^{(2+k)},\\
&&\nabla^pR_{\mathcal{M}_f}(X,Y,Y,X;Y,...,Y)R_{\mathcal{M}_f}(X,Y,Y,X)^{p-1}\nabla R_{\mathcal{M}_f}(X,Y,Y,X;Y)^{-p}\\
 &&\qquad=f^{(2+p)}\{f^{(2)}\}^{p-1}\{f^{(3)}\}^{-p}=\alpha_p(f)\,.
\end{eqnarray*}

Suppose $\phi:\mathcal{M}_{f_1}\rightarrow\mathcal{M}_{f_2}$ is a local isometry with $\phi(P_1)=P_2$. Set 
$$\Phi:=\phi_*(P_1):T_{P_1}M_1\rightarrow T_{P_2}M_2\,.$$
Assume that $f_1^{(3)}(P_1)\ne0$. Let
$X:=\Phi(\partial_x)$ and $Y:=(\Phi\partial_y)$. Since 
\begin{eqnarray*}
0&\ne&
R_{\mathcal{M}_{f_1}}(\partial_x,\partial_y,\partial_y,\partial_x;\partial_y)(P_1)=
(\Phi^*R_{\mathcal{M}_{f_2}})(\partial_x,\partial_y,\partial_y,\partial_x;\partial_y)(P_1)\\
&=&R_{\mathcal{M}_{f_2}}(X,Y,Y,X;Y)(P_2),
\end{eqnarray*}
we have $f^{(3)}(P_2)\ne0$. Let $p\ge2$. One may compute:
\begin{eqnarray*}\alpha_p(f_1,P_1)&=&\frac{\nabla^pR_{\mathcal{M}_{f_1}}(\partial_x,\partial_y,\partial_y,\partial_x;\partial_y,...,\partial_y)
    R_{\mathcal{M}_{f_1}}(\partial_x,\partial_y,\partial_y,\partial_x)^{p-1}}{
    \nabla R_{\mathcal{M}_{f_1}}(\partial_x,\partial_y,\partial_y,\partial_x;\partial_y)^{p}}(P_1)\\
&=&\frac{\nabla^pR_{\mathcal{M}_{f_2}}(X,Y,Y,X;Y,...,Y)R_{\mathcal{M}_{f_2}}(X,Y,Y,X)^{p-1}}{\nabla R_{\mathcal{M}_{f_2}}(X,Y,Y,X;Y)^{p}}(P_2)\\
&=&\alpha_p(f_2,P_2)\,.
\end{eqnarray*}
Thus Assertion (3a) implies Assertion (3b) in Theorem \ref{thm-1.6}.
Conversely, suppose 
$$\alpha_p(f_1,P_1)=\alpha_p(f_2,P_2)\quad\text{for}\quad p\ge2\,.$$
We use Lemma \ref{lem-3.2}
to see there is an isometry $\Phi$ from $T_{P_1}\mathcal{M}_{f_1}$ to $T_{P_2}\mathcal{M}_{f_2}$ so that
$$\Phi^*\nabla^kR_{\mathcal{M}_{f_2}}=\nabla^kR_{\mathcal{M}_{f_1}}\quad\text{for all}\quad k\,.$$
We may now apply
Lemma
\ref{lem-4.1} to see that $\mathcal{M}_{f_1}$ and $\mathcal{M}_{f_2}$ are isometric.
\end{proof}

To establish Theorem \ref{thm-1.4} (2), we shall need the following technical Lemma.
\begin{lemma}\label{lem-4.3} Let $f\in\mathcal{C}_3(\mathcal{O})$. If $\alpha_2(f)=k$ is constant, then either $f^{(2)}=ae^{\lambda
y}$ or 
$f^{(2)}=a(y+b)^c$.
\end{lemma}

\begin{proof} Let $h=f^{(2)}$. Then $h\ne0$, $h^\prime\ne0$,  and $h^{\prime\prime}h=kh^\prime h^\prime$. Thus
\begin{eqnarray*}
&&\textstyle\int\frac{h^{\prime\prime}}{h^{\prime}}=k\int\frac{h^\prime}h\quad\text{so}\quad
 \ln(h^\prime)=k\ln(h)+\beta\quad\text{so}\\
&&\textstyle h^\prime=e^\beta h^k\quad\text{so}\quad\int\frac{h^\prime}{h^k}=e^\beta y+\gamma\,.
\end{eqnarray*}
If $k=1$, this implies $\ln(h)=e^\beta y+\gamma$ or equivalently $h=e^\gamma e^{e^\beta y}$ which leads to an exponential solution.
If $k\ne1$, then
$h^{1-k}=(1-k)(e^\beta y+\gamma)$; this leads to a solution involving powers of a translate of $y$.
\end{proof}

\begin{proof}[Proof of Theorem \ref{thm-1.4} (2) and Theorem \ref{thm-1.8}] Suppose $f^{(2)}(y)=ae^{\lambda y}$. If $a=0$ or if
$\lambda=0$, then $\mathcal{M}_{f,P}$ is independent of $P$ and hence $\mathcal{M}_f$ is homogeneous by Lemma \ref{lem-4.1}. Thus we
may assume $a\ne0$ and $\lambda\ne0$ and hence $f\in\mathcal{C}_3$. Since $\alpha_p(f)$ is constant for $p\ge2$, Lemma \ref{lem-4.1}
implies
$\mathcal{M}_f$ is homogeneous.

Conversely, suppose $\mathcal{M}_f$ is homogeneous. If $\nabla R_{\mathcal{M}_f}=0$, then $f^{(2)}=c$ and we may take $\lambda=0$.
Thus we may assume $\nabla R_{{\mathcal{M}_f}}\ne0$ and hence $f^{(3)}\ne0$. Since the sign of $f^{(2)}$ determines the sign of
$R(\partial_x,\partial_y,\partial_y,\partial_y)$, we may suppose $f^{(2)}>0$; the case $f^{(2)}<0$ is entirely analogous. By Theorem
\ref{thm-1.6},
$\alpha_2$ is an isometry invariant. Thus $\alpha_2=k$ so $f^{(2)}=ae^{\lambda y}$ or $f^{(2)}=a(y+b)^c$. This latter case is ruled
out since it is not vanishing on $\mathbb{R}$; Theorem \ref{thm-1.4} (2) now follows.

The proof of Theorem
\ref{thm-1.8} follows the same lines but the additional solutions to the equation
$\alpha_2=k$ given by $a(y+b)^c$ now play a role.
\end{proof}

\section{Killing vector fields}\label{sect-5} We return to the generalized plane wave manifolds of Definition \ref{defn-2.1}. Let
$$\mathcal{K}_{\Xi,P}:=\{\eta\in T_PM:\nabla^kR_{\mathcal{P}_\Xi}(\eta,\xi_1,\xi_2,\xi_3;\xi_4,...,\xi_{k+3})=0\ 
\forall\ \xi_i\in T_P\mathbb{R}^{2p},\ \forall k\}\,.
$$
We have the following structure theorem:
\begin{theorem}\label{thm-5.1}  Suppose that $\mathcal{K}_{\Xi,P}=\operatorname{Span}\{\partial_i^{\tilde x}\}$ for all
$P\in\mathbb{R}^{2p}$. If
$X$ is a Killing vector field on
$\mathcal{M}_\Xi$, then there exists $\xi\in\mathbb{R}^p$, $A\in M_{p\times p}(\mathbb{R})$ and  a smooth
function
$\tilde\xi:\mathbb{R}^p\rightarrow\mathbb{R}^p$ so
that 
$$X(x,\tilde x)=\textstyle\sum_i(\xi_i+\textstyle\sum_jA_{ij}x_j)\partial_i^x
  -\textstyle\sum_i(\tilde\xi_i(x)+\textstyle\sum_jA_{ji}\tilde x_j)\partial_i^{\tilde x}\,.
$$
\end{theorem}

\begin{proof} We apply the discussion of Section \ref{sect-2} concerning the geodesics in $\mathcal{P}_\Xi$. Let $\phi$ be an
isometry of
$\mathcal{P}_\Xi$. Decompose
$\phi(x,0)=(\phi_1(x),\phi_2(x))$. Let 
\begin{eqnarray*}
&&\Phi(x):=\phi_*(x,0):T_{(x,0)}(\mathbb{R}^{2p})\rightarrow T_{(\phi_1(x),\phi_2(x))}(\mathbb{R}^{2p}),\\
&&\Phi\operatorname{Span}\{\partial_i^{\tilde x}\}=\Phi\mathcal{K}_{\Xi,(x,0)}
=\mathcal{K}_{\Xi,\phi(x,0)}=\operatorname{Span}\{\partial_i^{\tilde x}\}\,.
\end{eqnarray*}

Since $\gamma(t):=(x,t\tilde x)$ is a geodesic with $\gamma(0)=(x,0)$ and $\gamma^\prime(0)=(0,\tilde x)$ and
$\gamma_1:=\phi\circ\gamma$ is a geodesic with $\gamma_1(0)=(\phi_1(x),\phi_2(x))$ and
$\gamma_1^\prime(0)=(0,\Phi_1(x)\tilde x)$,
$$\phi(x,t\tilde x)=(\phi_1(x),\phi_2(x)+t\Phi_1(x)\tilde x)\,.$$

Fix $x$. Choose $\tilde x(t)$ so $\tau(t)=(tx,\tilde x(t))$
is a geodesic with $\tau(0)=(0,0)$ and $\tau^\prime(0)=(x,0)$. Thus $\phi\circ\tau$ is a geodesic starting at
$(\phi_1(0),\phi_2(0))$ with initial direction $(\Phi_2(0)x,\Phi_3(0)x)$ for suitably chosen $\Phi_2,\Phi_3\in M_p(\mathbb{R})$.
Thus
$$\phi\circ\tau(s)=(\phi_1(0)+s\Phi_2(0)x,\xi(s))$$ for some suitably chosen $\xi(s)$. Setting $t=1$ shows
$\phi_1(x)=\phi_1(0)+\Phi_2(0)x$. Thus
\begin{equation}\label{eqn-5.x}
\phi(x,\tilde x)=(\phi_1(x),\phi_2(x)+\Phi_1(x)\tilde x)=(\phi_1(0)+\Phi_2(0)x,\phi_2(x)+\Phi_1(x)\tilde x)\,.
\end{equation}

Let $\phi_s$ is a smooth $1$-parameter family of isometries. Differentiating Equation (\ref{eqn-5.x}) with respect
to $s$ and setting $s=0$ shows that Killing vectors have the form:
\begin{equation}\label{eqn-5.a}
X=\textstyle\sum_i(\xi_i+\sum_jA_{ij}x_j)\partial_i^x-\sum_i(\tilde\xi_i(x)+\sum_j\tilde
A_{ij}(x)\tilde x_j)\partial_i^{\tilde x}\,.\end{equation}

The Killing equation is
$$g_\Xi(\nabla_{\partial_i^x}X,\partial_j^{\tilde x})+
g_\Xi(\nabla_{\partial_j^{\tilde x}}X,\partial_i^x)=0$$
This equation then yields the relation $A_{ji}-\tilde A_{ij}(x)=0$. \end{proof}

\begin{remark}\label{rmk-5.2}\rm In deriving Equation (\ref{eqn-5.a}), we only needed the fact that $\phi$ was geodesic preserving
which is implied by the somewhat weaker assumption that
$\phi$ was an affine morphism, i.e. $\phi^*\nabla=\nabla$. We see therefore as a scholium that any affine Killing vector field on
$\mathcal{P}_\Xi$ has the form given in Equation (\ref{eqn-5.a}).
\end{remark}

We specialize Theorem \ref{thm-5.1} to the setting at hand:

\begin{theorem}\label{thm-5.3} If $f$ is not linear, then $X$ is a Killing vector field on $\mathcal{M}_f$ if and
only if
\begin{eqnarray*}
X&=&(\xi_1+A_{11}x+A_{12}y)\partial_x+(\xi_2+A_{21}x+A_{22}y)\partial_y\\
&-&(\tilde\xi_1(x,y)+A_{11}\tilde x+A_{21}\tilde y)\partial_{\tilde x}-(\tilde\xi_2(x,y)+A_{12}\tilde x+A_{22}\tilde y)
   \partial_{\tilde y}\end{eqnarray*}
where
\begin{eqnarray*}
0&=&-2fA_{11}-\partial_yf\cdot(\xi_2+A_{21}x+A_{22}y)-\partial_x\tilde\xi_1,\\
0&=&-2fA_{12}-\partial_x\tilde\xi_2-\partial_y\tilde\xi_1,\quad\text{and}\\
0&=&-\partial_y\tilde\xi_2\,.\end{eqnarray*}
\end{theorem}
\begin{proof} Let $X=\alpha\partial_x+\beta\partial_y-\tilde\alpha\partial_{\tilde x}-\tilde\beta\partial_{\tilde  y}$. If
$f$ is not linear, then there exists $k=k(P)$ so that $f^{(k)}(P)\ne0$ for some $k\ge2$. The non-degeneracy condition
$K_{\Xi,P}=\operatorname{Span}\{\partial_{\tilde x},\partial_{\tilde y}\}$ is satisfied as
$$\nabla^kR_{\mathcal{M}_f}(\partial_x,\partial_y,\partial_y,\partial_x;\partial_y,...,\partial_y)(P)\ne0\,.$$

The calculations of Section \ref{sect-2} show that
$$\begin{array}{lrrrr}
\nabla_{\partial_x}X=&(\partial_x\alpha)\partial_x&+(\partial_x\beta)\partial_y&
     +(-\partial_yf\cdot\beta-\partial_x\tilde\alpha)\partial_{\tilde
x}&+(\partial_yf\cdot\alpha-\partial_x\tilde\beta)\partial_{\tilde y}\\
\nabla_{\partial_y}X=&(\partial_y\alpha)\partial_x&+(\partial_y\beta)\partial_y&
     +(-\partial_yf\cdot\alpha-\partial_y\tilde\alpha)\partial_{\tilde x}&-(\partial_y\tilde\beta)\partial_{\tilde y}\\
\nabla_{\partial_{\tilde x}}X=&(\partial_{\tilde x}\alpha)\partial_x&+(\partial_{\tilde x}\beta)\partial_y&
     -(\partial_{\tilde x}\tilde\alpha)\partial_{\tilde x}&-(\partial_{\tilde x}\tilde\beta)\partial_{\tilde y}\\
\nabla_{\partial_{\tilde y}}X=&(\partial_{\tilde y}\alpha)\partial_x&+(\partial_{\tilde  y}\beta)\partial_y&
     -(\partial_{\tilde y}\tilde\alpha)\partial_{\tilde x}&-(\partial_{\tilde y}\tilde\beta)\partial_{\tilde y}\,.
\end{array}$$
As $X$ is a Killing vector field if and only if $g(\nabla_\xi X,\eta)+g(\nabla_\eta X,\xi)=0$ for all $\{\xi,\eta\}$,
$$\begin{array}{ll}
0=-2f\partial_x\alpha-\partial_yf\cdot\beta-\partial_x\tilde\alpha,&
0=-2f\partial_y\alpha-\partial_x\tilde\beta-\partial_y\tilde\alpha,\\
0=-2f\partial_{\tilde x}\alpha+\partial_x\alpha-\partial_{\tilde x}\tilde\alpha,&
0=-2f\partial_{\tilde y}\alpha+\partial_x\beta-\partial_{\tilde y }\tilde\alpha,\\
0=-\partial_y\tilde\beta,&
0=\partial_y\alpha-\partial_{\tilde x}\tilde\beta,\\
0=\partial_y\beta-\partial_{\tilde y}\tilde\beta,&
0=\partial_{\tilde x}\alpha,\\
0=\partial_{\tilde x}\beta+\partial_{\tilde y}\alpha,&
0=\partial_{\tilde y }\beta\,.
\end{array}$$
In view of
Theorem \ref{thm-5.1}, we may set
$$\begin{array}{ll}
\alpha=\xi_1+A_{11}x+A_{12}y,&\beta=\xi_2+A_{21}x+A_{22}y,\\
\tilde\alpha=\tilde\xi_1(x,y)+A_{11}\tilde x+A_{21}\tilde y,
&\tilde\beta=\tilde\xi_2(x,y)+A_{12}\tilde x+A_{22}\tilde y\,.\vphantom{\vrule height 12pt}
\end{array}$$
The desired results now follow. \end{proof}

\section{Killing vector fields on the manifolds $\mathcal{M}_f$}\label{sect-6}
 Let $G_{f,P}\subset G_f$ be the isotropy subgroup of isometries fixing a point $P$. Let 
$\gg_{f,P}$ and $\gg_f$ be the associated Lie algebras. Let $X$ be a Killing vector field on
$\mathcal{M}_f$ with
$X(P)=P$. Then $X$ generates a flow $\phi^t$ fixing $P$ and $d\phi^t(P)=e^{tA}$ generates a $1$-parameter family
of symmetries of the model
$\mathcal{U}_{\mathcal{M}_f,P}^\infty$ where $A\in\gg_{f,P}$ satisfies:
\begin{eqnarray*}
&&A:\partial_{x}\rightarrow A_{11}\partial_{x}+A_{21}\partial_{y}
-\partial_{x}\tilde\xi_1(P)\partial_{\tilde x}-\partial_{x}\tilde\xi_2(P)\partial_{\tilde y},\\
&&A:\partial_{y}\rightarrow A_{12}\partial_{x}+A_{22}\partial_{y}
-\partial_{y}\tilde\xi_1(P)\partial_{\tilde x}-\partial_{y}\tilde\xi_2(P)\partial_{\tilde y},\\
&&A:\partial_{\tilde x}\rightarrow-A_{11}\partial_{\tilde x}-A_{12}\partial_{\tilde y},\quad\text{and}\\
&&A:\partial_{\tilde y}\rightarrow-A_{21}\partial_{\tilde x}-A_{22}\partial_{\tilde y}\,.
\end{eqnarray*}

\begin{proof}[Proof of Theorem \ref{thm-1.11} (1)] Suppose that $f^{(2)}=0$. \label{sect-6.1}Then $R=0$ and hence $\mathcal{M}_f$ is
isometric to
$\mathbb{R}^{(2,2)}$. The isometry group is isometric to the warped product $O(2,2)\times\mathbb{R}^4$ and is $10$ dimensional as
claimed.\end{proof}

\begin{proof}[Proof of Theorem \ref{thm-1.11} (5)] Choose  a primitive $F$ so that $\partial_yF=f$ and $F(0)=0$. We exhibit theKilling
vector fields $X_i$, the associated flows $\Phi_t^{X_i}$, and the element of $A_i\in\gg_{f,0}$ and the symmetry $S_i:=e^{\varepsilon
A_i}$ for those vector fields that vanish at $0$.
 \begin{enumerate}

\item $X_1=\partial_x$;\quad ($\xi_1=1$);\quad
$\Phi_t^{X_1}:(x,y,\tilde x,\tilde y)\rightarrow(x+t,y,\tilde x,\tilde y)$.
\smallbreak\item $X_2=\partial_{\tilde x}$;\quad ($\tilde\xi_1=-1$);\quad
$\Phi_t^{X_2}:(x,y,\tilde x,\tilde y)\rightarrow(x,y,\tilde x+t,\tilde y)$.
\smallbreak\item $X_3=\partial_{\tilde y}$;\quad ($\tilde\xi_2=-1$);\quad
$\Phi_t^{X_3}:(x,y,\tilde x,\tilde y)\rightarrow(x,y,\tilde x,\tilde y+t)$.
\smallbreak\item $X_4=-y\partial_{\tilde x}+x\partial_{\tilde y}$;\quad  ($\tilde\xi_1=y$, $\tilde\xi_2=-x$);\hfill\break
$\Phi_t^{X_4}:(x,y,\tilde x,\tilde y)\rightarrow(x,y,\tilde x-ty,\tilde y+tx)$;\hfill\break
$A_4:\partial_x\rightarrow\partial_{\tilde y}$,\quad $A_4:\partial_y\rightarrow-\partial_{\tilde x}$,\quad
$A_4:\partial_{\tilde x}\rightarrow 0$,
$A_4:\partial_{\tilde y}\rightarrow 0$,\hfill\break
$S_4:\partial_x\rightarrow\partial_x+\varepsilon\partial_{\tilde y},$\quad
$S_4:\partial_y\rightarrow\partial_y-\varepsilon\partial_{\tilde x},$\quad
$S_4:\partial_{\tilde x}\rightarrow\partial_{\tilde x},$\quad
$S_4:\partial_{\tilde y}\rightarrow\partial_{\tilde y}$.
\smallbreak\item $X_5=y\partial_x+2F\partial_{\tilde x}-\tilde x\partial_{\tilde y}$;\quad
($A_{12}=1$, $\tilde\xi_1=-2F$);\hfill\break
$\Phi_t^{X_5}:(x,y,\tilde x,\tilde  y)\rightarrow(x+ty,y,\tilde x+2F(y)t,\tilde y+t\tilde x+F(y)t^2)$.\hfill\break 
$A_5:\partial_x\rightarrow0$,\quad
$A_5:\partial_y\rightarrow\partial_x-f\partial_{\tilde x}$,\quad
$A_5:\partial_{\tilde x}\rightarrow-\partial_{\tilde y}$,\quad
$A_5:\partial_{\tilde y}\rightarrow0$,\hfill\break
$S_5:\partial_x\rightarrow\partial_x+\varepsilon f\partial_{\tilde y}$,\quad
$S_5:\partial_y\rightarrow\partial_y+\varepsilon\partial_x+2\varepsilon f\partial_{\tilde x}-\varepsilon^2f\partial_{\tilde y}$,
\hfill\break
$S_5:\partial_{\tilde x}\rightarrow\partial_{\tilde x}-\varepsilon\partial_{\tilde y}$,\quad
$S_5:\partial_{\tilde y}\rightarrow\partial_{\tilde  y}$.
\end{enumerate}

Suppose that $f^{(2)}$ is non-constant. Thus $f^{(3)}$ must be non-vanishing at some point $P$; to simplify the notation, we may
assume without loss of generality that $P=0$. Since
$\mathcal{U}_{f,0}$ is isomorphic to $\mathcal{U}^1$, $\dim\gg_{f,P}\le2$ by Lemma \ref{lem-3.3}. Since $X_4$ and $X_5$ are Killing
vector fields vanishing at $0$ with $A_4$ and $A_5$ linearly independent,
$\dim\gg_{f,P}=2$. If $f\ne ae^{\lambda y}$ for $a\ne0$ and if $f\ne a(y+b)^c$ for $a\ne0$ and $c\ne 0,1,2$, then $\alpha_1$ is
non-constant near $0$ and hence $\mathcal{M}_f$ is not locally homogeneous at $0$. $\dim\gg_f(0)\le3$ and the $5$ Killing
vector fields listed above are a basis for $\gg_f$.
\end{proof}

\begin{proof}[Proof of Theorem \ref{thm-1.11} (3,4)] Also suppose that $f=ae^{\lambda y}$ or $f=a(y+b)^c$. As we can choose $P$ with
$f^{(3)}(P)\ne0$, $\mathcal{U}_{f,P}^1$ is isomorphic to $\mathcal{U}^1$ and thus $\dim\gg_{f,P}\le2$ by Lemma \ref{lem-3.3}. Thus
the argument given above shows $\dim\gg_f\le6$ and to complete the proof, we must only exhibit an additional vector field.
\begin{enumerate}
\item Suppose that $f(y)=e^{\lambda y}$. Set\hfill\break
$X_6:=-\frac\lambda2x\partial_x+\partial_y+\frac\lambda2\tilde x\partial_{\tilde x}$;\quad
($A_{11}=-\frac\lambda2$,$\xi_2=1$).\hfill\break
$\Phi_t^{X_6}:(x,y,\tilde x,\tilde y)\rightarrow(e^{-\frac\lambda2t}x,y+t,e^{\frac\lambda2t}\tilde x,\tilde y)$.
\item Suppose $f=a(y+b)^c$. By renormalizing our coordinates, we may suppose $a=1$ and $b=0$. Set\hfill\break
$X_6:=x\partial_x-\frac2cy\partial_y-\tilde x\partial_{\tilde x}+\frac2c\tilde
y\partial_{\tilde y}$;\quad $(A_{11}=1,A_{22}=-\frac2c)$;\hfill\break
$\Phi_t^{X_6}(x,y,\tilde x,\tilde y):=(e^{t}x,e^{-\frac2ct}y,e^{-t}\tilde x,e^{{\frac2c}t}\tilde y)$.
\end{enumerate}
This completes the proof.
\end{proof}
\begin{remark}\label{rmk-6.1}\rm If $f=ae^{\lambda y}$, then the flows defined by $X_1$, $X_2$, $X_3$, and $X_6$ act transitively 
on $\mathbb{R}^4$. This gives a direct proof that $\mathcal{M}_{f}$ is a homogeneous space. If $f=ay^n$ for $n\in\mathbb{N}$, then
the flow defined by $X_6$ fixes the hyperplane $y=0$. Let 
$$\mathcal{O}:=\{(x,y,\tilde x,\tilde y):y>0\}\,.$$
The flows defined by $X_1$, $X_2$,
$X_3$, and
$X_6$ define a transitive action on
$\mathcal{M}_{f,\mathcal{O}}$; thus $\mathcal{M}_{f,\mathcal{O}}$ is a homogeneous proper open incomplete submanifold of
$\mathcal{M}_f$. Such examples can not exist in the Riemannian setting.
\end{remark}

\begin{proof}[Proof of Theorem \ref{thm-1.11} (2)] By rescaling, we may suppose $f(y)=\pm y^2$; we suppose $f=y^2$ as the case
$f=-y^2$ is similar. We then have
$\mathcal{U}_{f,0}^0$ is isomorphic to
$\mathcal{U}^0$ and thus by Lemma \ref{lem-3.3}, $\dim\gg_{f,0}\le4$. Consequently $\dim\gg_f\le8$. To establish the desired
result, we must construct 8 additional Killing vector fields.
\begin{enumerate}
\item $X_6:=\partial_y+2xy\partial_{\tilde x}-x^2\partial_{\tilde y}$;\quad $(\xi_2=1$, $\tilde\xi_1=-2xy$,
$\tilde\xi_2=x^2$);\hfill\break
$\Phi_t^{X_6}:(x,y,\tilde x,\tilde y)\rightarrow(x,y+t,\tilde x+ 2xyt+xt^2,\tilde y-x^2t)$.
\smallbreak\item $X_7:=x\partial_y+(yx^2-\tilde y )\partial_{\tilde x}-\frac13x^3\partial_{\tilde y}$;\quad
($A_{21}=1$, $\tilde\xi_1=-yx^2$, $\tilde\xi_2=\frac13x^3$);\hfill\break
$A_7:\partial_x\rightarrow\partial_y$,\quad $A_7:\partial_y\rightarrow 0$,\quad
$A_7:\partial_{\tilde x}\rightarrow0$,\quad $A_7:\partial_{\tilde y}\rightarrow-\partial_{\tilde x}$;\hfill\break
$S_7:\partial_y\rightarrow\partial_y$,\quad
$S_7:\partial_{\tilde x}\rightarrow\partial_{\tilde x}$,\quad
$S_7:\partial_{\tilde y}\rightarrow\partial_{\tilde y}-\varepsilon\partial_{\tilde x}$.
\smallbreak\item $X_8:=x\partial_x-y\partial_y-\tilde x\partial_{\tilde x}+\tilde y\partial_{\tilde y}$;\quad $(A_{11}=1$,
$A_{22}=-1$);
\hfill\break
$A_8:\partial_x\rightarrow\partial_x$,\quad $A_8:\partial_y\rightarrow-\partial_y$,\quad
$A_8:\partial_{\tilde x}\rightarrow-\partial_{\tilde x}$, $A_8:\partial_{\tilde y}\rightarrow\partial_{\tilde  y}$;
\hfill\break
$S_8:\partial_x\rightarrow e^{\varepsilon}\partial_x$,\quad
$S_8:\partial_y\rightarrow e^{-\varepsilon}\partial_y$,\quad
$S_8:\partial_{\tilde x}\rightarrow e^{-\varepsilon}\partial_{\tilde x}$,\quad
$S_8:\partial_{\tilde y}\rightarrow e^{\varepsilon}\partial_{\tilde y}$.
\end{enumerate}
The proof is now complete.
\end{proof}

\begin{remark}\label{rmk-6.2}\rm The vector field $X_6$ generates the missing translational symmetry; the flows
for $X_1,X_2,X_3,X_6$ act transitively on $\mathbb{R}^4$; this gives a direct proof that $\mathcal{M}_{y^2}$ is a homogeneous space.
Furthermore, the isomorphisms generated by $X_4,X_5,X_7,X_8$ generate the full symmetry group of the model $\mathcal{U}^0$.
We have omitted the flows for $X_7$ and $X_8$ in the interests of brevity.
\end{remark}

\section*{Acknowledgements} Research of Peter Gilkey partially supported by the Atlantic Association for Research in the
Mathematical Sciences (Canada) and by the Max Planck Institute in the Mathematical Sciences (Leipzig, Germany). 
Research of S. Nik\v cevi\'c partially supported by MM 1646 (Serbia).

\end{document}